 \documentclass[11pt]{amsart}
\usepackage{geometry}   
\usepackage{setspace}
\usepackage{graphicx}
\usepackage{amssymb, amsthm}
\usepackage{epstopdf}
\usepackage{pinlabel}
\usepackage{color}
 \newtheorem{theorem}{Theorem}[section]

    \newtheorem{proposition}[theorem]{Proposition}

    \theoremstyle{definition}
\newtheorem{definition}[theorem]{Definition}

\newtheorem{example}[theorem]{Example}
\newtheorem{examples}[theorem]{Examples}

  \newcommand{\Z}{\ensuremath{{\mathbb{Z}}}}
\newcommand{\R}{\ensuremath{{\mathbb{R}}}}
\newcommand{\Q}{\ensuremath{{\mathbb{Q}}}}
  
\newcommand{\E}{\ensuremath{{\mathbb{E}}}}
\newcommand{\Hy}{\ensuremath{{\mathbb{H}}}}
  
      \newcommand{\T}{\ensuremath{{\mathbb{T}}}}

 \newcommand{\mbX}{\ensuremath{{\partial_M X}}}

\newcommand{\mbY}{\ensuremath{{\partial_M Y}}}
\newcommand{\mbG}{\ensuremath{{\partial_M G}}}

\title{Searching for Hyperbolicity}
\author{Ruth Charney}

\thanks {R. Charney was partially supported by NSF grant DMS-1607616}

\begin{document}

\begin{abstract} 
This paper is an expanded version of a talk given at the AWM Research Symposium 2017.  It is intended as a gentle introduction to geometric group theory with a focus on the notion of hyperbolicity, a theme that has inspired the field from its inception to current-day research.
\end{abstract}

\maketitle

\section{Introduction}
This paper is an expanded version of a talk given at the AWM Research Symposium 2017.  It is intended as a gentle introduction to geometric group theory for the non-expert, with a focus on the notion of hyperbolicity.  Geometric group theory came into its own in the 1990's, in large part due to a seminal paper of Mikhail Gromov \cite{Gro}.   While the field has grown considerably since that time, hyperbolicity remains a central theme and continues to drive much current research in the field. 

As the name suggests, geometric group theory provides a bridge between groups viewed as algebraic objects and geometry.  Group arise in all areas of mathematics and can be described in many different ways.  Some arise purely algebraically (such as certain matrix groups), others have combinatorial descriptions (via presentations), and still others are defined topologically (such as fundamental groups of topological spaces).  Geometric group theory is based on the principle that if a group acts as symmetries  of some geometric object, then one can use geometry to better understand the group.  

For many groups, it is easy to find such an action.  The symmetric group on $n$-letters acts by symmetries on an $n$-simplex and the dihedral group of order $2n$ is the symmetry group of an $n$-gon.  This is also the case for some infinite groups. The free abelian group $\Z^n$ acts by translation on $\R^n$ (preserving the Euclidean metric), and the free group on $n$-generators act by translation on a regular tree of valence $2n$ with edges of length one.  

So the first question one might ask is, when can one find such an action?  Given an abstract group $G$, can we always realize $G$ as a group of symmetries of some geometric object?  As we will see below, the answer is yes.  However, some geometric objects are more useful than others for this purpose.  In the early 1900's, Max Dehn was interested in groups arising as fundamental groups of hyperbolic surfaces.  These groups act by isometries on the hyperbolic plane $\Hy^2$.  Dehn used the geometry of the hyperbolic plane to prove some amazing properties of these groups. We will discuss one of these results in Section \ref{sec:hyperbolic} below.  Decades later, these ideas motivated Gromov, who introduced the notion of a hyperbolic metric space.  He showed that the properties that Dehn deduced held more generally for any group acting nicely on such a space.

Gromov's notion of hyperbolic spaces and hyperbolic groups have been much studied since that time.  Many well-known groups, such as mapping class groups and fundamental groups of surfaces with cusps, do not meet Gromov's criteria, but nonetheless display some hyperbolic behavior.  In recent years, there has been much interest in capturing and using this hyperbolic behavior wherever and however it occurs.  

In this paper, I will review some basic notions in geometric group theory, discuss Dehn's work and Gromov's notion of hyperboicity, then introduce the reader to some recent developments in the search for hyperboliclity.  My goal is to be comprehensible, not comprehensive.  For those interested in learning more about the subject, I recommend \cite{Bo06} and \cite{ClMa} for a general introduction to geometric group theory and \cite{CDP} and  \cite {GhHa} for more about hyperbolic groups.

\section{Geodesic metric spaces, isometries and quasi-isometries}

We begin with some basic definitions.
Let $X$ be a metric space with distance function $d: X \times X \to \R$.   
A \emph{geodesic} in $X$ is a distance preserving map from an interval $I \subset \R$  into $X$, that is, a map $\alpha : I \to X$ such that for all $t_1,t_2 \in I$, 
$$d(\alpha(t_1),\alpha(t_2)) = |t_i -t_2|. $$  
The interval $I$ may be finite or infinite.  This definition is analogous to the notion of a geodesic in a Riemannian manifold.  In particular, a geodesic between two points in $X$ is a length-minimizing path.

A \emph{geodesic metric space} is a metric space $X$ in which any two points are connected by a geodesic.  For such metric spaces, the distance is intrinsic to the space; the distance between any two points is equal to the minimal length of a path connecting them.  Often, we also require that our metric space be \emph{proper}, that is, closed balls in $X$ are compact.

\begin{examples}  (1)  Consider the unit circle $S^1$ in the plane.  There are two natural metrics we could put on $S^1$.  The first is the induced Euclidean metric: the distance between two points is the length of the straight line in $\R^2$ between them.  The other is the arc length metric:  the distance between two points is the length of the (shortest) circular arc between them.  The first of these is not a geodesic metric (since, for example, there is no path in $S^1$ of length 2 connecting a pair of antipodal points) whereas the second one is geodesic.  

(2) Suppose $\Gamma$ is a connected graph.  There is a natural geodesic metric on $\Gamma$ obtained by identifying each edge with a copy of the unit interval $[0,1]$ and defining the distance between any two points in $\Gamma$ to be the length of the shortest path between them.  This metric is proper if and only if each vertex has finite valence.

(3) Let $M$ be a complete Riemannian manifold.  Then the usual distance function given by minimizing path lengths is a proper, geodesic metric on $M$.  
\end{examples}

A map between two metric spaces $f : X \to Y$ is an \emph{isometry} if it is bijective and preserves distances.  In lay terms, an isometry of $X$ to itself is a ``symmetry" of $X$.  These symmetries form a group under composition.

Now suppose we are given a group $G$.  Our first goal is to find a nice metric space on which $G$ acts as a group of symmetries.  
Sometimes, such an action arises naturally.  For example, suppose $G$ is the fundamental group of a
Riemannian manifold M.  Then passing to the universal cover $\widetilde M$, we get an action of $G$ by deck transformations on $\widetilde M$.  This action is distance preserving since it takes geodesic paths to geodesic paths. 

More generally, the same works for the fundamental group of any geodesic metric space $X$ that admits a universal cover.  The universal cover $\widetilde X$ inherits a geodesic metric such that the projection to $X$ is a local isometry, and the deck transformations act isometrically on $\widetilde X$.

\begin{example}  Consider the free group on two generators $F_2$.  This group is the fundamental group of a wedge of two circles, $S^1 \vee S^1$, so it acts by isometries on the universal cover, namely the regular 4-valent tree.  
\end{example}

In general, a group $G$ can act by isometries on a variety of different geodesic metric spaces.  Some of these actions, however, are not helpful in studying the group.  For example, any group acts on a single point!  To have any hope that the geometry of the space will produce information about the group, we will need some extra conditions on the action.

\begin{definition}  A group $G$ is said to act \emph{geometrically} on a metric space $X$ if the action satisfies the following three properties.
\begin{itemize}
\item \emph{isometric:}  Each $g \in G$ acts as an isometry on $X$.
\item \emph{proper:}  For all $x \in X$ there exits $r > 0$ such that 
$\{g \in G \mid B(x,r) \cap gB(x,r) \neq \emptyset \}$  is finite, where $B(x,r)$ denotes the ball of radus $r$ centered at $x$. 
\item \emph{cocompact:}  There exists a compact set $K \subset X$ whose translates by $G$ cover $X$
(or equivalently, $X/G$ is compact).
\end{itemize}
\end{definition}

In particular, the fundamental group of a compact metric space $X$ acts geometrically on the universal covering space $\widetilde X$. 
But now suppose that our group $G$ arises purely algebraically.  How can we find a metric space $X$ on which $G$ acts geometrically?  Here is a construction that works for any finitely generated group. 
 
Choose a finite generating set $S$ for $G$.  Define the \emph{Cayley graph} for $G$ with respect to $S$ to be the graph $\Gamma_S(G)$ whose vertices are in one-to-one correspondence with the elements of $G$ and for each $s \in S, g \in G,$ there is an edge (labelled by $s$) connecting the vertex $g$ to the vertex $gs$.  

$G$ acts on $\Gamma_S(G)$ by left multiplication on the vertices.  That is, $h \in G$ maps the vertex labelled $g$ to the vertex labelled $hg$.  Note that $h$ takes edges to edges, the edge connecting $g$ to $gs$ maps to the edge connecting $hg$ to $hgs$.  Thus, if we put the path-length metric on $\Gamma_S(G)$ as described above, then this action preserves the metric and is easily seen to be geometric.  

\begin{figure}[htpb] 
\centering
\includegraphics[height=4 cm]{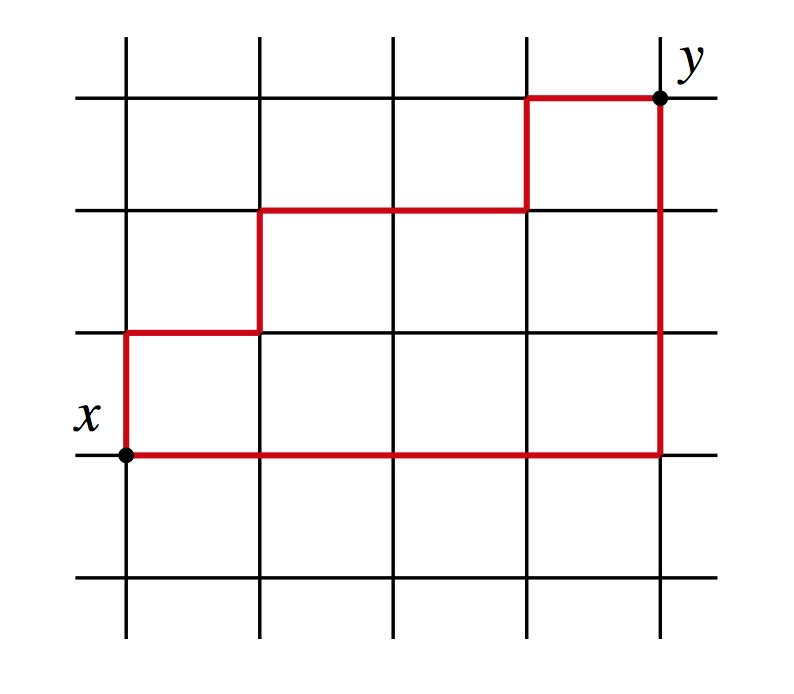}
\caption{Geodesics from $x$ to $y$ in the Cayley graph of $\Z^2$ }\label{fig:grid}
\end{figure}

\begin{figure}[htpb] 
\centering
\includegraphics[height=4 cm]{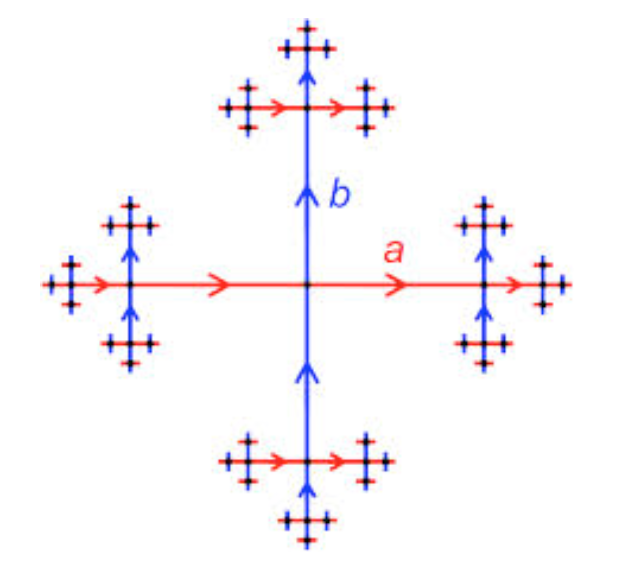}
\caption{Cayley graph of $F_2$ }\label{fig:F2-tree}
\end{figure}

\begin{example}  (1) Consider the free abelian group $\Z^2$ with generating set $S=\{(1,0), (0,1)\}$.  Viewing $\Z^2$ as the set of integer points in the plane $\R^2$, we can identify the Cayley graph of $\Z^2$ with the square grid connecting these points.  Distances  are measured by path lengths in this grid, so the distance from $(0,0)$ to $(n,m)$ is $|n|+|m|$.  Note that there are, in general, many geodesic paths between any two points.  See Figure \ref{fig:grid}.

(2) Let $F_2$ be the free group with generating set $S=\{a,b\}$.  View $F_2$ as the fundamental group of the wedge of two circles labelled $a$ and $b$.  Lifting these labels to the universal covering space, we get a 4-valent tree, as in Figure \ref{fig:F2-tree}, with every red edge labelled $a$ and every blue edge labelled $b$  (Metrically, you should picture every edge as having the same length.)
 This tree is precisely the Cayley graph $\Gamma_S(F_2)$.  To see this, choose a base vertex $v$ and identify each vertex with the element of $F_2$ that translates $v$ to that vertex.  
\end{example}


Clearly, the Cayley graph depends on the choice of generating set.  For example, were we to add a third generator $(1,1)$ to our generating set for $\Z^2$, the Cayley graph would get additional edges which cross the grid diagonally and would hence change the distance between vertices.  

This could be a cause for concern; we seek geometric properties of the Cayley graph that are intrinsic to the group, so they should not be dependent on a choice of generating set.  Luckily,
in the case of a finitely generated group, replacing one finite generating set by another, does not distort  distances too badly.  This leads to a fundamental concept in geometric group theory.

\begin{definition}  A map $f : X \to Y$ between two metric spaces is a \emph{quasi-isometric embedding} if there exists constants $K,C$ such that for all $x,z \in X$
$$\frac{1}{K}\,  d_X(x,z) - C \leq  d_Y(f(x),f(z)) \leq K\, d_X(x,z) + C.$$
If in addition, every point in $Y$ lies within $C$ of some point in $f(X)$, then $f$ is a  \emph{quasi-isometry}. In this case we write $X \sim_{QI} Y$.
\end{definition}

It can be shown that any quasi-isometry has a ``quasi-inverse", so the relation $X \sim_{QI} Y$ is an equivalence relation.  We remark that quasi-isometries need not be continuous maps.

\begin{figure}[htpb] 
\centering
\includegraphics[height=3.5 cm]{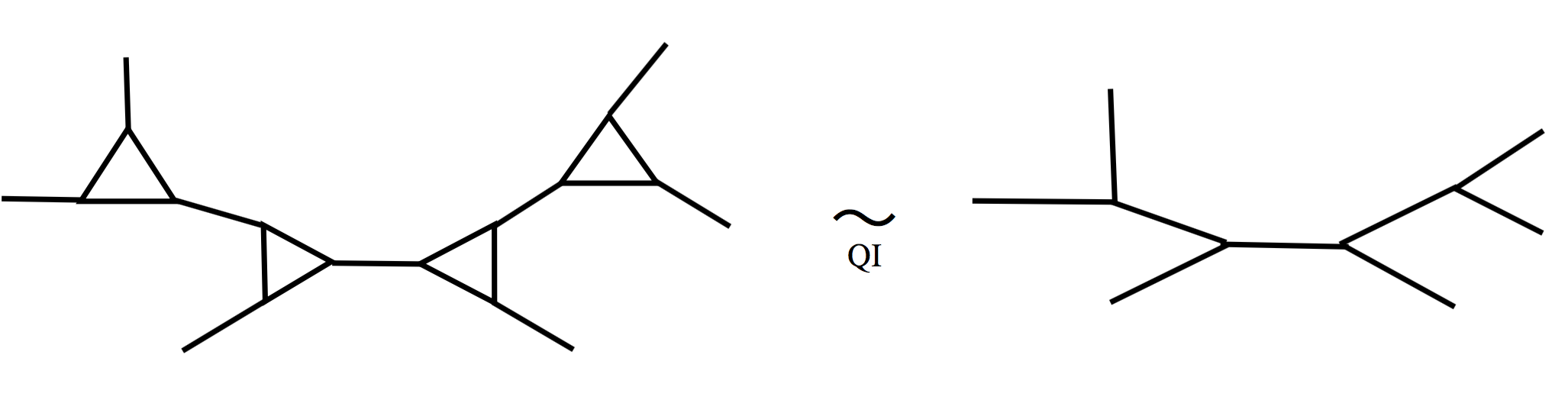}
\caption{A graph quasi-isometric to a tree}\label{fig:QI}
\end{figure}

\begin{examples} (1) Consider the inclusion of the integer grid into the plane $\R^2$.  This is a quasi-isometry.  The quasi-inverse is a discontinuous map sending the interior of each square to its boundary.

(2)  Consider the graph in Figure \ref{fig:QI}.  Collapsing each of the triangles to a point gives a quasi-isometry of this graph onto the 3-valent tree.  We call such a graph is a quasi-tree.
\end{examples}

It is an easy exercise to show that if $G$ is a finitely generated group, then the Cayley graphs with respect to different finite generating sets are all quasi-isometric.  (In fact, they are bi-Lipschitz, i.e., we can take the maps to be continuous and the constant $C$ to be zero.)  Thus, identifying $G$ with the vertex set in its Cayley graph, we can view $G$ itself as a metric space and this metric is well-defined up to quasi-isometry.  
In fact, we have the following more general statement.

\begin{proposition}[\v{S}varc-Milnor Lemma] Suppose $G$ acts geometrically on a geodesic metric space $X$.  Then $G$ is finitely generated and for any choice of basepoint $x_0 \in X$, the map $G \to X$ taking $g \mapsto gx_0$ is a quasi-isometry.
\end{proposition}

As a result, the notion of quasi-isometry is a fundamental concept in geometric group theory.  The properties  of a group that one can hope to glean from its action on a metric space are generally properties preserved by quasi-isometries.  Moreover, since any finite index subgroup $H < G$ is quasi-isometric to $G$, often the strongest statement we can make regarding a property $(P)$ is that our group $G$ is \emph{virtually} (P), meaning that some finite index subgroup of $G$ satisfies property $(P)$.  

The classification of finitely generated groups up to quasi-isometry is a meta-problem in the field.  It is easy to see that any two groups that are commensurable (i.e., they contain subgroups of finite index that are isomorphic) are quasi-isometric.  So a related problem is the question of rigidity: for a given group $G$, is every group quasi-isometric to $G$ also commensurable to $G$?   

\section{Hyperbolic groups}\label{sec:hyperbolic}

Once we have our group acting geometrically on a metric space, we can ask how geometric properties of the space are reflected in algebraic or combinatorial properties of the group. The classical example of this comes from the work of Max Dehn \cite{De}.  Dehn was interested in fundamental groups of surfaces. A closed orientable surface of genus $g \geq 2$ (i.e. a torus with $g$ holes) can be given a Riemannian metric of constant curvature $-1$ and its universal covering space can be identified with the hyperbolic plane $\Hy^2$.  This gives a geometric action of the fundamental group on $\Hy^2$.
Using geometric properties of the hyperbolic plane, Dehn proved a some very strong combinatorial properties for these groups.  I will describe two of his results here.  

One way to describe a group is by means of a presentation.  Given a set of generators $S$ for a group $G$, there is a natural surjection from the free group $F(S)$ onto $G$.  The kernel $K$ of this map is a normal subgroup of $F(S)$.  Formally, a \emph{presentation} of $G$ consists of a generating set $S$, together with a set $R \subset F(S)$ such that $R$ generates $K$ as a normal subgroup.  In other words, $G \cong F(S) / \langle \langle R \rangle\rangle$.  We denote such a presentation for $G$ by writing
$$G = \langle S \mid R \rangle.$$

In practice, we usually indicate $R$ by a set of equations that hold in $G$.  Viewing elements of $F(S)$ as ``words" in the alphabet $S \cup S^{-1}$, the elements in $R$ are words that are equal to the identity in $G$.  However, many other words, such as products and conjugates of those in $R$, are also equal to the identity in $G$.  The idea is that \emph{all} relations among the generators that hold in $G$ should be consequences of the ones listed in the presentation. 

Here are some examples. The cyclic group of order $n$ has presentation
$$ \Z / n\Z = \langle s \mid s^n=1 \rangle $$
while the free abelian group on two generators has presentation
$$ \Z^2 =\langle a,b \mid ab=ba \rangle. $$
Can you recognize the following group?
$$ G = \langle u,v \mid u^4=1, u^2=v^3 \rangle.$$
It turns out that this group is isomorphic to the special linear group $SL(2,\Z)$!  The isomorphism is given by identifying
$$u = \begin{pmatrix}
0 & -1\\
1 & 0
\end{pmatrix} \quad\quad
v=\begin{pmatrix}
0 & -1\\
1 & 1
\end{pmatrix}.
$$


Every group can be described by a presentation, though in general $S$ and $R$ need not be finite.
Presentations can be extremely useful, and are the starting point for combinatorial group theory.  On the other hand, presentations can sometimes be very mysterious and frustratingly difficult to decifer. For example, consider the following questions.  
\begin{enumerate}
\item \emph{The Word Problem: } Given a finite presentation $\langle S \mid R \rangle$ and a word $w$ in $F(S)$, is there an algorithm to decide whether $w$ represents the trivial element in $G$?
\item \emph{The Isomorphism Problem:}  Given two finite presentations $\langle S \mid R \rangle$ and $\langle S' \mid R' \rangle$, is there an algorithm to decide whether the groups they represent are isomorphic?
\end{enumerate}
It turns out that there are groups for which no such algorithms exist!  In this case we say the Word Problem (or the Isomorphism Problem) is unsolvable.  What Dehn showed, using the geometry of hyperbolic space, was that for fundamental groups of hyperbolic surfaces, both of these problems are solvable.  Moreover, he showed that the word problem is solvable in linear time; he described an algorithm (now known as a Dehn's algorithm) for which the time it takes to decide whether  $w \in F(S)$ represents the identity element in $G$, is linear in the length of $w$  viewed as a word in $S \cup S^{-1}$.

Some 75 year later, Mikhail Gromov made a startling observation:  the only property of the hyperbolic plane that Dehn really needed to derive his results was the fact that triangles in $\Hy^2$, no matter how far apart their vertices may be, are always ``thin".  And from this observation, the modern field of hyperbolic geometry (and more generally geometric group theory) was born.  

What do we mean by ``thin"? 
Let $X$ be any geodesic metric space.  A triangle $T(a,b,c)$ in $X$ consists of three vertices $a,b,c \in X$ together with a choice of geodesics connecting them.  

\begin{figure}[htpb] 
\centering
\includegraphics[height=5 cm]{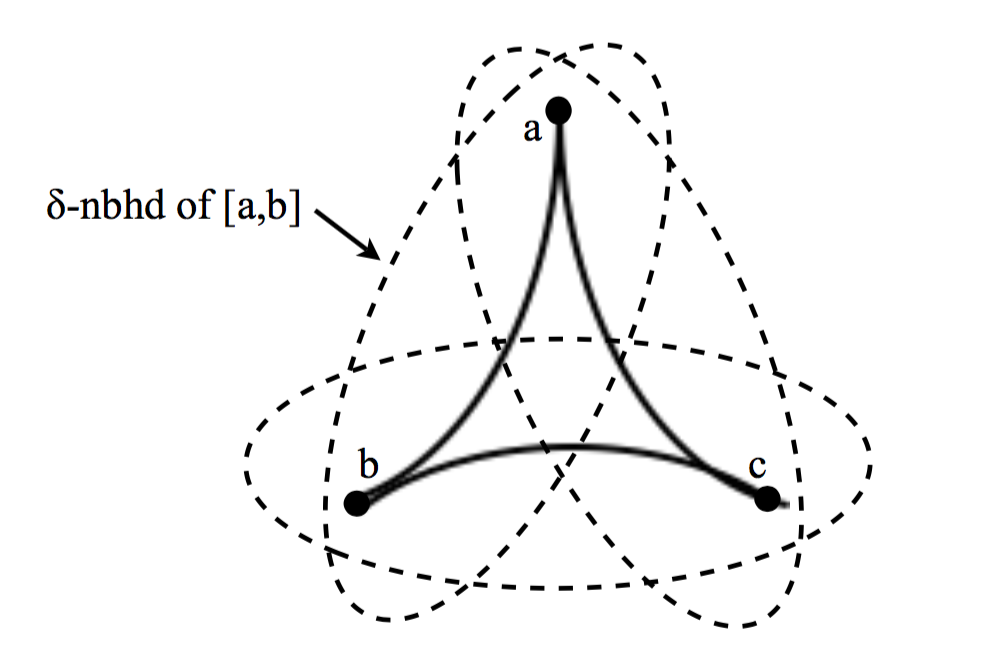}
\caption{A $\delta$-thin triangle}\label{fig:hyp triangle}
\end{figure}

\begin{definition}  Let $X$ be a geodesic metric space and let $\delta \geq 0$.  We say a triangle $T(a,b,c)$ in $X$ is \emph{$\delta$-thin} if each side of $T$ lies in the $\delta$-neighborhood of the other two sides. (See Figure \ref{fig:hyp triangle}.)
\end{definition}

One can show that in $\Hy^2$, every triangle, even those with vertices at infinity (ideal triangles) are $\delta$-thin for $\delta = \ln(1 + \sqrt{2})$.  This fact was crucial to Dehn's work.

This brings us finally to Gomov's notion of hyperbolicity.  

\begin{definition}  A geodesic metric space $X$ is \emph{$\delta$-hyperbolic} if every triangle in $X$ 
is $\delta$-thin.  We say $X$ is \emph{hyperbolic} if it is $\delta$-hyperbolic for some $\delta$.  A finitely generated group $G$ is \emph{hyperbolic} if it acts geometrically on a hyperbolic metric space.
\end{definition}

One can show that if two spaces are quasi-isometric and one of them is hyperbolic, then so is the other (though the constant $\delta$ may change).  In particular, a finitely generated group $G$ is hyperbolic if and only if some (hence any) Cayley graph of $G$ is hyperbolic.  

\begin{figure}[htpb] 
\centering
\includegraphics[height=3.5 cm]{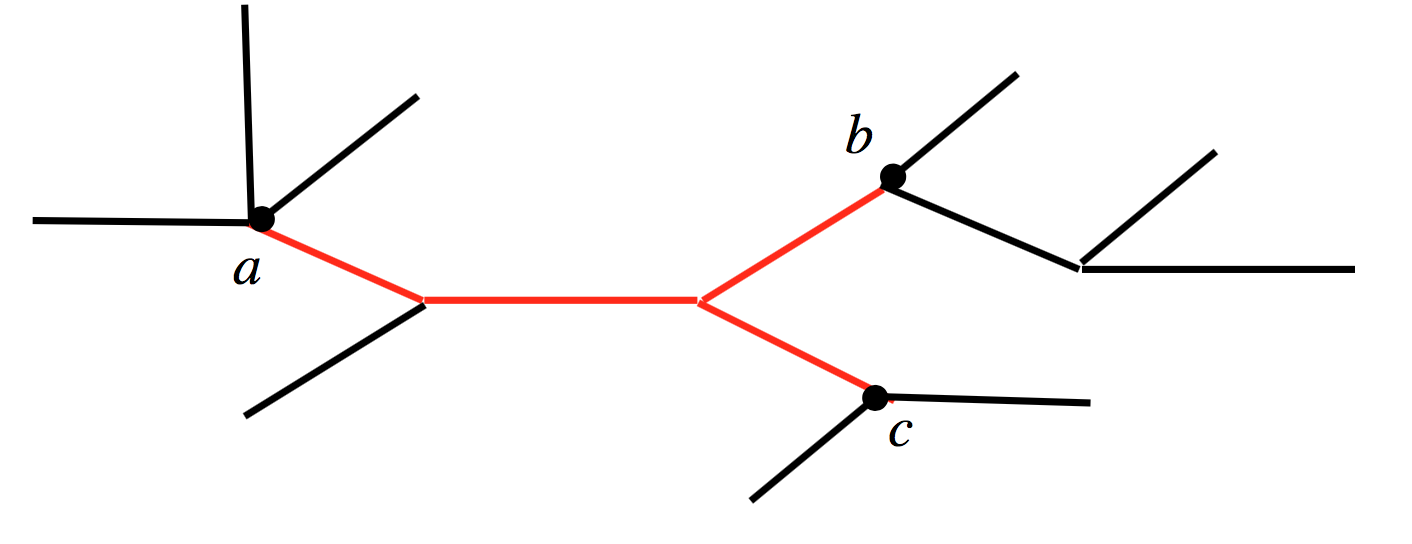}
\caption{A triangle in a tree is 0-thin.}\label{fig:tree}
\end{figure}

\begin{examples}
(1) First, a trivial example.  Any bounded metric space $X$ is $\delta$-hyperbolic where $\delta$ is the diameter of $X$, and hence any finite group is hyperbolic.  

(2) Let $X$ be an infinite tree.  Then for any three points $a,b,c$ in $X$, the triangle connecting them degenerates into a tripod and is hence 0-hyperbolic!  (See figure  \ref{fig:tree}.)   Since the Cayley graph of a free group is a tree, it follows that finitely generated free groups are hyperbolic.

(3)  Recall the presentation of $SL(2,\Z)$ given above.  The center of $SL(2,\Z)$ is the order two subgroup generated by $u^2 = v^3$.  Modding out by this subgroup gives the group $PSL(2,\Z)$ with presentation 
$$ PSL(2,\Z) = \langle u,v \mid u^2=v^3=1 \rangle.$$
The Cayley graph of $PSL(2,\Z)$ with respect to this generating set is the quasi-tree drawn in Figure \ref{fig:QI} (continued out to infinity), with the edges of  triangles labelled $v$ and the remaining edges labeled $u$.  $SL(2,\Z)$ also acts geometrically on this quasi-tree (with the center acting trivially),  $SL(2,\Z)$ and $PSL(2,\Z)$ are both hyperbolic groups.   

(4)  Here is a non-example.  Let $\R^2$ be the plane with the standard Euclidean metric.  Taking larger and larger isosceles right triangles, we can see that there is no bound on ``thinness".  Since the Cayley graph of $\Z^2$ is quasi-isometric to $\R^2$, $\Z^2$ is not hyperbolic.  Indeed, it is a theorem that a hyperbolic group cannot contain a copy of  $\Z^2$.
\end{examples}

Now suppose that we are given a hyperbolic group $G$.  What does the geometry tell us about the group?  Here is a list of some consequences of hyperbolicity.  We refer the reader to \cite{BrHa} and \cite{CDP} for proofs and additional references. 
\begin{enumerate}
\item  $G$ has a finite presentation.
\smallskip
\item  $G$ has a Dehn's algorithm, hence a linear time solution to the Word Problem.
\smallskip
\item The Isomorphism Problem is solvable for the class of hyperbolic groups.
\smallskip
\item  The centralizer of every element of $G$ is virtually cyclic.
\smallskip
\item $G$ has at most finitely many conjugacy classes of torsion elements.
\smallskip
\item For any finite set of elements $g_1, \dots g_k$ in $G$, there exists $n>0$ such that the set 
 $\{g_1^n \dots g_k^n\} $ generates a free subgroup of rank at most $k$ in $G$.
\smallskip
\item  For $n$ sufficiently large, $H^n(G;\Q) =0$.
\smallskip
\item  If $G$ is torsion-free, it has a finite $K(G,1)$-space. 
\end{enumerate}

The proofs of these properties are  beyond the scope of this paper, but the conclusion should be clear:  geometry can have strong implications for algebraic and combinatorial properties of a group!

\section{Beyond hyperbolicity}

Classical hyperbolic geometry and Gromov's generalization to $\delta$-hyperbolic spaces  have provided powerful tools for studying hyperbolic groups.  But this class of groups is very special.  For example, any group containing a subgroup isomorphic to $\Z^2$ cannot be hyperbolic (see property (4) above).   In recent years, there has been much interest in generalizing some of these techniques to broader classes of groups.  Gromov himself introduced a notion of ``non-positive curvature" for geodesic metric spaces, called CAT(0) spaces, and groups acting on these spaces have also been extensively studied.  CAT(0) geometry, for example, played a major role in the recent work of Agol, Wise, and others leading to a proof of the Virtual Haken Conjecture, the last remaining piece in Thurston's program to classify three-manifolds.  But this is a topic for another day.

Other approaches to generalizing the theory of hyperbolic groups involve looking at groups that act on hyperbolic spaces, but where the actions are not geometric; instead they satisfy some weaker conditions.  This includes, for example, ``relatively hyperbolic groups",  ``acylindrically hyperbolic groups", and ``hierarchically hyperbolic groups".   The first of these is modeled on fundamental groups of hyperbolic manifolds with cusps;  the others are inspired by mapping class groups and their actions on curve complexes.  Some very nice introductions to these topics can be found in \cite{Bo12, Osi, Sis}.

My own work in this area has focused on a somewhat different approach to capturing hyperbolic behavior in more general spaces and groups.  Let's begin with an example.  

\begin{figure}[htpb] 
\centering
\includegraphics[height=7 cm]{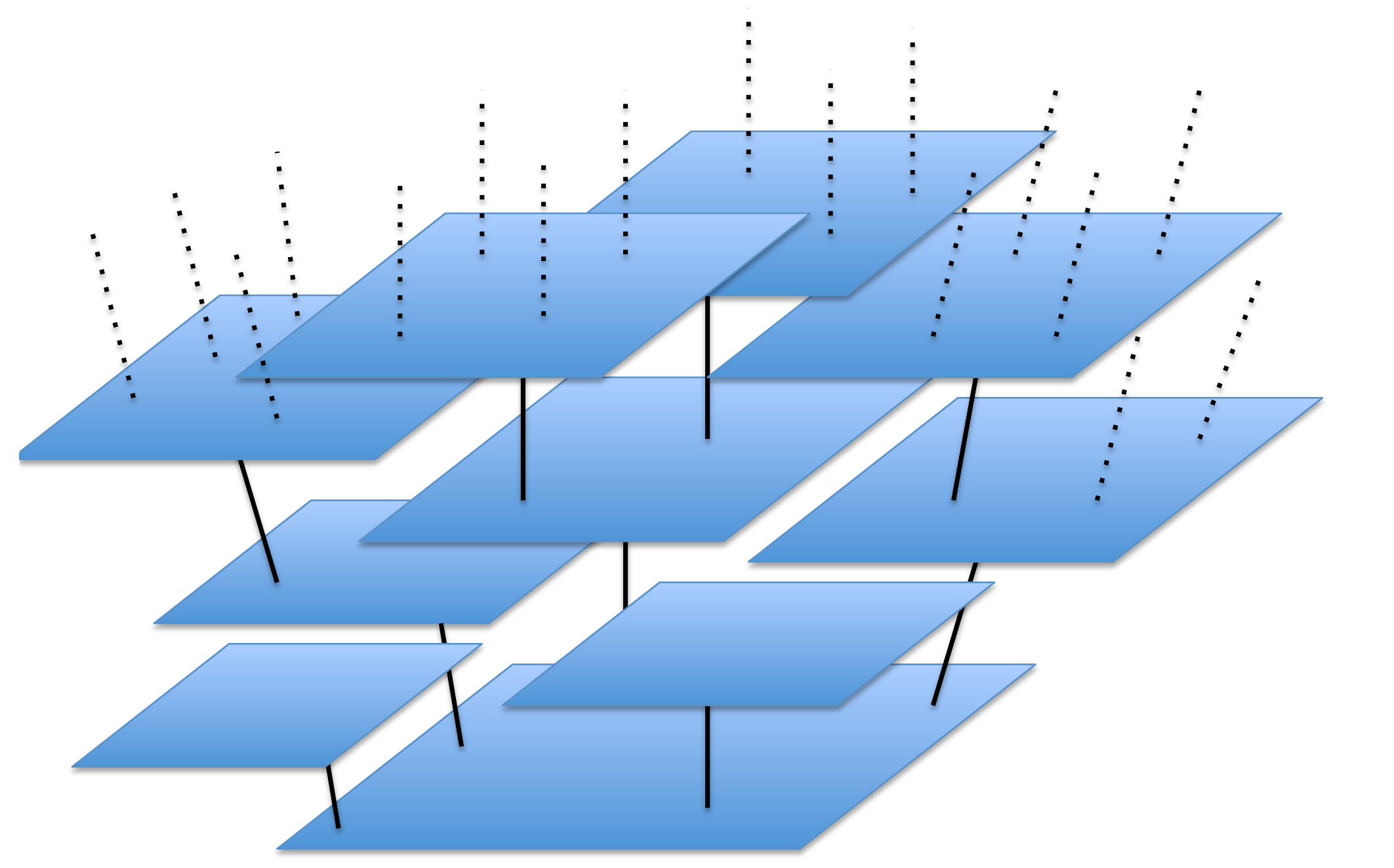}
\caption{Tree of flats}\label{fig:flats}
\end{figure}

\begin{example}\label{tree of flats}  
Let $Z$ be the space obtained by gluing a circle and a torus together at a single point,
$$Z = S^1 \vee T^2.$$
Let $\widetilde Z$ be its universal cover.  In $\widetilde Z$, covering the torus we see infinitely many copies of the Euclidean plane (which we refer to as ``flats") and emanating from each lattice point in each of these planes is a line segment covering the circle.  The result is a tree-like configuration of planes and lines (see Figure \ref{fig:flats}), which we will refer to as the ``tree of flats".  Triangles lying in a single flat can be arbitrarily ``fat" whereas triangles whose sides lie mostly along the vertical lines are ``thin".  Thus, in $\widetilde Z$, we have hyperbolic-like directions, and non-hyperbolic directions.  Intuitively, the more we travel vertically, the more hyperbolic it feels.
\end{example}

We are interested in identifying geodesics in a metric space $X$ that behave like geodesics in a hyperbolic space.  A good way to encode such geodesics is by means of a boundary.  
In general, unbounded metric spaces do not come equipped with a boundary.  For example, in the hyperbolic plane (or the Euclidean plane) one can travel forever in any direction. To create a boundary for such a space, we need to add a point for each ``direction to infinity".  In the case of the hyperbolic or Euclidean plane, this space of directions forms a circle.
Adding this circle to the plane ``at infinity" compactifies the space. 

It turns out that this idea generalizes nicely to any hyperbolic metric space.
For a $\delta$-hyperbolic space $X$, we define the boundary as follows.  
A \emph{ray} in $X$ is an isometric embedding 
$\alpha: [0,\infty) \to X$.  As a set, the boundary of $X$ is defined to be
$$\partial X =\{ \alpha \mid \alpha: [0,\infty) \to X \,\,\textrm{is a ray} \} / \sim $$
where $\alpha \sim \beta$ if  $\alpha$ and $\beta$ remain bounded distance from each other.  
In the Euclidean plane, for example, two rays are equivalent if and only if they are parallel.

To topologize $\partial X$, think of two rays as representing nearby points in the boundary if they remain close to each other for a long time.  More precisely, define a neighborhood $N(\alpha, R)$ of a ray $\alpha$ to be the set of rays $\beta$ such that $\beta(t)$ lies within $2\delta$ of $\alpha(t)$ for $0 \leq t \leq R$.  As $R$ increases, these neighborhood get smaller and smaller, and together, they form a neighborhood basis for a topology on $\partial X$.  

For example, the boundary of the hyperbolic plane $\Hy^2$ is a circle while the boundary of an infinite tree is a Cantor set.

From the geometric group theory viewpoint, a key property of the boundary of a hyperbolic space is quasi-isometry invariance.

\begin{theorem}  Let $f : X \to Y$ be a quasi-isometry between two hyperbolic metric spaces.  Then $f$ induces a homeomorphism $\partial f: \partial X \cong \partial Y$.  In particular, a hyperbolic group $G$ has a well-defined boundary, namely the boundary of a Cayley graph of $G$.  
\end{theorem}

These boundaries have many nice properties and applications.  The boundary gives rise to a compactification of $X$, $\overline X = X \cup \partial X$, and it provides a powerful tool for studying the dynamics of groups actions, rigidity theorems, geodesic flows, etc.  

In the quest to extend the techniques of hyperbolic geometry to more general spaces and groups, it is natural to ask whether analogous boundaries can be defined in more general contexts.  Certainly we can consider equivalence classes of geodesic rays in any geodesic metric space $X$.  However, if $X$ is not hyperbolic, many things can go wrong.  In some cases, it is not even clear how to define a topology on this set as the neighborhoods described above need not satisfy the requirements for a neighborhood basis.  

Moreover, even when there is a nice topology on $\partial X$, other fundamental properties of hyperbolic boundaries can fail to hold. Consider, for example, the boundary of the Euclidean plane.  This boundary is a circle and  each point on the boundary can be represented by a unique ray based at the origin. As observed above, it provides a compactification of the plane. That's the good news.  Here is some bad news.
\begin{itemize}
\item  Many isometries (in particular all translations) act trivially on the boundary.  
\item  The only pairs of points on the boundary that can be joined by a bi-infinite geodesic are pairs of antipodal points.  
\item  A quasi-isometry of the plane to itself need not extend to a map on the boundary.  For example,
the map $f: \R^2 \to \R^2$ taking $re^{i \theta} \mapsto re^{i(\theta + \ln(r))}$ is a quasi-isometry which twists each ray emanating from the origin into a spiral.
\end{itemize}
In short, many of the properties of hyperbolic boundaries that permit applications to dynamics, rigidity, etc., fail to hold for this boundary.

Most significantly from the point of view of geometric group theory, quasi-isometry invariance fails miserably for non-hyperbolic boundaries.  There are examples of groups that act geometrically on two CAT(0) spaces (spaces of non-positive curvature) whose boundaries are not homeomorphic \cite{CrKl}.  Thus, we don't have a well-defined notion of a boundary for these groups.

What goes wrong is the failure of the Morse property.  A quasi-isometry $f : X \to Y$ of hyperbolic spaces takes a geodesic ray in $X$ to a quasi-geodesic ray in $Y$, that is, a quasi-isometric embedding of the half-line $\R^+=[0,\infty)$ into $Y$.  The Morse property guarantees that this quasi-geodesic ray lies close to some geodesic ray and hence determines a well-defined point at infinity.  

\begin{definition} A ray (or bi-infinite geodesic) $\alpha$ in $X$ is \emph{Morse} if there exists a function $N: \R^+ \times \R^+ \to \R^+$ such that for any $(K,C)$-quasi-geodesic $\beta$ with endpoints on $\alpha$, $\beta$ lies in the $N(K,C)$-neighborhood of $\alpha$.  The function $N$ is called a \emph{Morse gauge} for $\alpha$ and we say that $\alpha$ is $N$-Morse.
\end{definition}

If $X$ is hyperbolic, then there exists a Morse gauge $N$ such that every ray in $X$ is $N$-Morse. 
This property is the key to proving quasi-isometry invariance for boundaries of hyperbolic spaces,
and it plays a key role in the proofs of many other properties of hyperbolic spaces as well. If $X$ is not hyperbolic, the Morse property will fail for some (or perhaps all) rays in $X$.  On the other hand, many non-hyperbolic spaces contain a large number of Morse rays.  Consider, for example, our ``tree of flats" described above.  It can be shown that a ray that spends a uniformly bounded amount of time in any flat is Morse.  We view Morse rays as ``hyperbolic-like" directions in $X$.  Indeed, it can be shown that these rays share many other nice properties with rays in hyperbolic space \cite{ChSu, CoHu}.  For example, if two sides of a triangle are $N$-Morse, then the triangle is $\delta$-thin where $\delta$ depends only on the Morse gauge $N$. 

Which brings us to the Morse boundary.  The \emph{Morse boundary}, $\mbX$ can be defined for any proper geodesic metric space $X$.  As a set, it consists of equivalence classes of Morse rays,
$$ \mbX = \{ \alpha \mid \textrm{$\alpha : [0,\infty) \to X$ is a Morse ray}\} / \sim$$
where the equivalence $\sim$ is defined as before.  The topology is more subtle.  For a sequence of rays $\{\alpha_i\}$ to converge to $\alpha$ in this topology, they must not only converge pointwise, they must also be uniformly Morse, that is, there exists a Morse gauge $N$ such that all of the $\alpha_i$ are $N$-Morse. 


This boundary was first introduced for CAT(0) spaces by myself and Harold Sultan in \cite{ChSu} and shown to be quasi-isometry invariant.  This was then generalized to arbitrary proper geodesic metric spaces by Matt Cordes in \cite{Co15}. 

\begin{theorem}  Let $f: X \to Y$ be a quasi-isometry between two proper geodesic metric spaces.  Then $f$ induces a homeomorphism $\partial f: \mbX \to \mbY$.  In particular, $\mbG$ is well-defined for any finitely generated group. 
\end{theorem}

\begin{examples}
(1)  If $X$ is hyperbolic, then all rays are $N$-Morse for some fixed $N$.  So in this case,
$\mbX = \partial X$, is the usual hyperbolic boundary.  For example,  
$\partial_M \Hy^2$ is a circle. 

(2)  For $X= \E^2$ the Euclidean plane, there are no Morse rays at all, so $\mbX= \emptyset$

(3)  Let $\widetilde Z$ be the ``tree of flats" from Example \ref{tree of flats}.  Then the Morse geodesics in $\widetilde Z$  are those that spend a uniformly bounded amount of time in any flat and the maximum time spent in a flat is determined by the Morse gauge. Thus, for a given Morse gauge $N$, the $N$-Morse rays emanating from a fixed basepoint $z_0$, lie in a subspace quasi-isometric to a  tree.  It follows that these rays determine a Cantor set in the boundary and the Morse boundary of $\widetilde Z$  is the direct limit of these Cantor sets.  Note that since $\widetilde Z$ is the universal cover of $S^1 \vee T^2$, the fundamental group $\Z * \Z^2 = \pi_1(S^1 \vee T^2)$  acts geometrically on $\widetilde Z$.  Hence by the theorem above, the Morse boundary of the group $\Z * \Z^2$ is homeomorphic to the Morse boundary of $\widetilde Z$.
\end{examples}

The Morse boundary was designed to capture hyperbolic-like behavior in non-hyperbolic metric spaces
and to give a well-defined notion of a boundary for a finitely generated group $G$.  When the Morse boundary is non-trivial, it provides a new tool for studying these spaces and groups.  It can be used, for example, to study the dynamics of isometries \cite{Mur} and to determine when two groups are quasi-isometric \cite{ChMu}.  It can also be used to study geometric properties of subgroups $H < G$ \cite{CoHu}.  For a survey of recent results on Morse boundaries see \cite{Co17}.  

\bigskip
Geometric group theory is a broad and growing area of mathematics.  This article is intended only as a snapshot of some themes that run through the field.  I invite you to investigate further!

\bibliographystyle{siam}

\end{document}